\documentclass{amsart}
\usepackage{amsfonts, amsthm, amssymb, amsmath, stmaryrd}
\usepackage{mathrsfs,array}
\usepackage{eucal,color,times,enumerate,accents}
\usepackage[all]{xy}
\usepackage{url}
\usepackage{enumitem}
\usepackage{tikz}
\usepackage{tikz-cd}
\usepackage[pdftex,bookmarksnumbered,bookmarksopen]{hyperref}
\usepackage{stackrel}
\hypersetup{colorlinks, linkcolor=blue}
\input xy
\xyoption{all}
\newcommand{\Zl}{\mathbb{Z}_{\ell}}

\newcommand{\Ql}{\mathbb{Q}_{\ell}}

\newcommand{\Xbar}{X_{\overline{k}}}

\newcommand{\Gm}{\mathbb{G}_m}

\newcommand{\F}{\mathbb{F}}
\newcommand{\Q}{\mathbb{Q}}
\newcommand{\Z}{\mathbb{Z}}

\newcommand{\etale}{\'etale }

\usepackage{calligra}
\makeatletter
\def\@splitop#1#2\@nil{$\mathscr{#1}\!\!\!$\calligra#2\;\!}
\newcommand*\DeclareCursiveOperator[2]{%
  \newcommand#1{\mathop{\mbox{\@splitop#2\@nil}}\nolimits}}
\makeatother
\DeclareCursiveOperator{\EXT}{Ext}
\DeclareCursiveOperator{\HOM}{Hom}

\begin{document}
\bibliographystyle{alpha}
\newtheorem{theorem}[equation]{Theorem}

\newtheorem{fact}[equation]{Fact}

\newtheorem{proposition}[equation]{Proposition}
\newtheorem{lemma}[equation]{Lemma}
\newtheorem{corollary}[equation]{Corollary}
\newtheorem{claim}[equation]{Claim}
\newtheorem{question}[equation]{Question}
\newtheorem{conjecture}[equation]{Conjecture}
\newtheorem{counterexample}[equation]{Counterexample}

\newtheorem{answer}[equation]{Answer}
\newtheorem{example}[equation]{Example}
\newtheorem{warning}[equation]{Warning}
\newtheorem{notation}[equation]{Notation}
\newtheorem{construction}[equation]{Construction}

\theoremstyle{definition}
\newtheorem{remark}[equation]{Remark}
\newtheorem{definition}[equation]{Definition}

\newcommand{\adjunction}[4]{\xymatrix@1{#1{\ } \ar@<-0.3ex>[r]_{ {\scriptstyle #2}} & {\ } #3 \ar@<-0.3ex>[l]_{ {\scriptstyle #4}}}}

\title{Uniform boundedness for Brauer groups of forms in positive characteristic}
\author{Emiliano Ambrosi}
\numberwithin{equation}{subsubsection} 

\begin{abstract}
Let $k$ be a finitely generated field of characteristic $p>0$ and $X$ a smooth and proper scheme over $k$. Recent works of Cadoret, Hui and Tamagawa show that, if $X$ satisfies the $\ell$-adic Tate conjecture for divisors for every prime $\ell\neq p$, the Galois invariant subgroup  $Br(X_{\overline k})[p']^{\pi_1(k)}$ of the prime-to-$p$ torsion of the geometric Brauer group of $X$ is finite. The main result of this note is that, for every integer $d\geq 1$, there exists a constant $C:=C(X,d)$ such that for every finite field extension $k \subseteq k'$ with $[k':k]\leq d$ and every $(\overline k/k')$-form $Y$ of $X$ one has $|(Br(Y\times_{k'}\overline k)[p']^{\pi_1(k')}|\leq C$. The theorem is a consequence of general results on forms of compatible systems of $\pi_1(k)$-representations and it extends to positive characteristic a recent result of Orr and Skorobogatov in characteristic zero.
\end{abstract}
\maketitle

\tableofcontents
\numberwithin{equation}{subsubsection} 
\section{Introduction}
Let $k$ be a field of characteristic $p\geq 0$ with algebraic closure $\overline k$ and write $\pi_1(k)$ for the absolute Galois group of $k$. In this paper, a $k$-variety is a reduced scheme, separated and of finite type over $k$ and, if $X$ is a $k$-variety, we write $\Xbar:=X\times_k \overline k$. The letter $\ell$ will always denote a prime $\neq p$. 
\subsection{Brauer groups}
\subsubsection{Finiteness of Brauer groups}
Let $X$ be a $k$-variety. Write $Br(X_{\overline k})[p']$ for the prime-to-$p$ torsion of the (cohomological) Brauer group $Br(X_{\overline k}):=H^{2}(X_{\overline k},\Gm)$ of $X_{\overline k}$ and recall that if $X$ is smooth over $k$ then $Br(X_{\overline k})$ is a torsion group. If $k$ is finitely generated and $X$ is smooth and proper over $k$, one expects $Br(X_{\overline k})[p']^{\pi_1(k)} $ to be small. This is predicted by (variants of) the $\ell$-adic Tate conjecture for divisors (\cite{tateconj}):
\begin{conjecture}[$T(X,\ell)$]
Assume that $k$ is finitely generated and $X$ is a smooth and proper $k$-variety. Then the $\ell$-adic cycle class map
$$c_{X_{\overline k}}:Pic(X_{\overline k})\otimes \Ql\rightarrow \bigcup_{[k':k]<+\infty}H^2(X_{\overline k},\Ql(1))^{\pi_1(k')}$$
is surjective.
\end{conjecture}
As it is well known (see e.g. \cite[Proposition 2.1.1]{brauer}), Conjecture $T(X,\ell)$ holds if and only if, for any finite field extension $k\subseteq k'$, the $\ell$-primary torsion $Br(X_{\overline k})[\ell^{\infty}]^{\pi_1(k')}$ of $Br(X_{\overline k})^{\pi_1(k')}$ is finite. But one can expect stronger finiteness results.  
\begin{fact}\label{finiteness}
Assume that $k$ is finitely generated and $X$ is a smooth and proper $k$-variety. Then:
\begin{enumerate}
\item \cite[Theorem 5.5]{skoro}: If $p=0$ and the integral Mumford Tate conjecture for $X$ holds (\cite[Conjecture
C.3]{integraltate}), then $Br(X_{\overline k})^{\pi_1(k)}$ is finite;
\item \cite[Corollary 1.5]{ultra}: If $p>0$ and $T(X,\ell)$ holds for every prime $\ell\neq p$ (or equivalently for one prime $\ell\neq p$), then $Br(X_{\overline k})[p']^{\pi_1(k)}$ is finite.
\end{enumerate}
\end{fact}
\subsubsection{Uniform boundedness in forms}
Let $X$ be a smooth proper variety over a finitely generated field $k$. Recall that for a field extension $k\subseteq k'\subseteq \overline k$, a $(\overline k/k')$-form of $X$ is a $k'$-variety $Y$ such that $Y_{\overline k}:=Y\times_{k'}\overline k\simeq \Xbar$. Let $k\subseteq k'$ be a finite field extension and let $Y$ be a $(\overline k/k')$-form of $X$. If $p=0$ and $X$ satisfies the integral Mumford Tate conjecture (resp. if $p>0$ and $T(X,\ell)$ holds for every prime $\ell\neq p$), then the same is true for $Y$, hence $Br(Y_{\overline k})^{\pi_1(k)}$ (resp. $Br(Y_{\overline k})[p']^{\pi_1(k')}$) is a finite group. But, for an integer $d\geq 1$, it is not clear whether one can find a uniform bound (depending only on $X$ and $d$) for $|Br(Y_{\overline k})^{\pi_1(k')}|$ (resp. $|Br(Y_{\overline k})[p']^{\pi_1(k')}|$), while $k'$ is varying among the finite field extensions $k\subseteq k'$ with $[k':k]\leq d$ and $Y$ among the $(\overline k/k')$-forms of $X$. If $p=0$, this is proved by Orr-Skorobogatov in \cite[Theorem 5.1]{skoro}. If $p>0$, this is the first main result of this note.
\begin{theorem}\label{mainBrauer}
Assume that $k$ is finitely generated, $X$ is a smooth proper $k$-variety and $p>0$. If $T(X,\ell)$ holds for every prime $\ell\neq p$ (or equivalently for one prime $\ell\neq p$), then for every integer $d\geq 1$, there exists a constant $C:=C(X,d)$ such that for every finite field extension $k \subseteq k'$ of degree $\leq d$ and every $(\overline k/k')$-form $Y$ of $X$ one has
$$(Br(Y_{\overline k})[p'])^{\pi_1(k')}\leq C.$$
\end{theorem}
\subsection{Forms of representations}
\numberwithin{equation}{subsubsection} 
Theorem \ref{mainBrauer} is a consequence of two general results (Propositions \ref{maximalcompatible} and \ref{adeliccompatible}) on compatible systems of $\pi_1(k)$-representations. Before stating them, we introduce some definitions and notation. In the following, $k$ is a finitely generated field of characteristic $p>0$, $\F_q$ (resp. $\F$) is the algebraic closure of $\F_p$ in $k$ (resp. in $\overline k$) and we write $k_\F:=k\otimes_{\F_q}\F\simeq k\F\subseteq \overline k$. Set $\ell_{0}=3$ (resp. $\ell_0=2$) if $p\neq 3$ (resp. $p=3$) and $s_{\ell}=\ell$ (resp. $s_{\ell}=4$) if $\ell\neq 2$ (resp. $\ell=2$). Fix a collection $\underline T:=\{T_{\ell}\}_{\ell\neq p}$ of rank $r$ finitely generated $\Zl$-modules endowed with a continuous action of $\pi_1(k)$.
\subsubsection{Definitions}
We say that $\underline T$ is a compatible system of $\pi_1(k)$-modules if there exists a smooth geometrically connected $\F_q$-variety $\mathcal K$ with generic point $Spec(k)\rightarrow \mathcal K$ such that, for every prime $\ell\neq p$, the action of $\pi_1(k)$ on $T_{\ell}$ factors trough the canonical surjective morphism $\pi_1(k)\twoheadrightarrow \pi_1(\mathcal K)$ and the collection $\{{V_{\ell}:=T_{\ell}\otimes \Ql}\}_{\ell\neq p}$ give rise to a $\Q$-rational compatible system on $\mathcal K$ in the sense of Serre: for each closed point $\mathfrak t\in \mathcal K$, the characteristic polynomial of the arithmetic Frobenius at $\mathfrak t$ acting on $V_{\ell}$ is in $\Q[T]$ and independent of $\ell$.
\begin{remark}
The notion of compatible system is stable under subquotients and the usual operations $\oplus$, $\otimes$, $^\vee$.
\end{remark}
\begin{definition}\label{defcomp}
Let $k\subseteq k'$ be a finite field extension. A $(\overline k/k')$-form of $\underline T$ is a compatible system of $\pi_1(k')$-representations $\underline U$ such that, for each $\ell\neq p$, there exists a finite field extension $k'\subseteq k_{\ell}$ and an isomorphism of $\pi_1(k_{\ell})$-modules $T_{\ell}\simeq U_{\ell}$.
\end{definition}
\subsubsection{Results}
In Definition \ref{defcomp}, the extension $k\subseteq k_{\ell}$ is allowed to depend on $\ell$. Our first main result in this setting produces an extension of (explicitly) bounded degree that works for every prime $\ell\neq p$. Let $?\in \{\emptyset,\F\}$.
\begin{proposition}\label{maximalcompatible}
Let $\underline U$ be a $(\overline k/k)$-form of $\underline T$. Then, there exists a finite field extension $k_?\subseteq k_{\underline U}$ of degree $\leq |GL_{r}(\Z/s_{\ell_0})|^2$ and a $\pi_1(k_{\underline U})$-equivariant isomorphism $T_{\ell}/(T_{\ell})_{tors}\simeq U_{\ell}/(U_{\ell})_{tors}$ for every prime $\ell\neq p$.
\end{proposition}
Proposition \ref{maximalcompatible} reduces the problem of bounding uniformly the invariants of forms of $\underline{T}$ to studying the action of $\pi_1(k')$ on $\underline{T}$, when $k\subseteq k'$ is varying among the finite field extensions of bounded degree. In this setting we prove:
\begin{proposition}\label{adeliccompatible}
Suppose that $T_{\ell}$ is torsion free for $\ell\gg 0$. Then there exists a finite field extension $k_?\subseteq k'$ of degree $\leq |GL_{r}(\Z/s_{\ell_0})|$ with the following property:
For every integer $d\geq 1$ there exists a constant $C:=C(\underline T,d)$ such that, for every finite field extension $k'\subseteq k''$ of degree $\leq d$, one has 
$$\prod_{\ell\neq p}[(T_{\ell}\otimes \Ql/\Zl)^{\pi_1(k'')}:(T_{\ell}\otimes \Ql/\Zl)^{\pi_1(k')}]\leq C.$$
\end{proposition}
\begin{remark}\label{inutile}
In the proof of Theorem \ref{mainBrauer} we only use the version of Propositions \ref{maximalcompatible} and \ref{adeliccompatible} where $?=\emptyset$. On the other hand, the proofs of the two versions are very similar and we believe that both versions are of independent interest.
\end{remark}
\subsection{Motivic representation}
The main motivation to state Theorems \ref{maximalcompatible} and \ref{adeliccompatible} in this generality is that they apply directly to representations associated to $\ell$-adic \etale cohomology of smooth proper $k$-varieties; see Subsections \ref{sectionmotivicbrauer} and \ref{Sectionabelian}. Since Theorems \ref{maximalcompatible} and \ref{adeliccompatible} require only the compatibility of the compatible system and not further assumptions as purity, one could apply them also to representations arising from the cohomology of some not necessarily smooth and proper $k$-varieties (e.g. semi-abelian schemes). 
\subsection{Strategy}
To prove Proposition \ref{maximalcompatible}, first we prove a group theoretic proposition (Proposition \ref{connectedcomponent}) that bounds the number of connected components of the Zariski closure of the image of an $\ell$-adic representation of a profinite group, only in terms of $\ell$ and of the rank of the representation. To get Proposition \ref{maximalcompatible}, one has to get rid of the dependency on $\ell$. This follows formally from the fact that the connectedness of the $\ell$-adic monodromy group can be read off the L-function of the various compatible systems $\{T_{\ell}^{\otimes n}\otimes (T^{\vee}_{\ell})^{\otimes m}\}_{\ell\neq p}$.

For the proof of Proposition \ref{adeliccompatible}, the key point is to show that, if the Zariski closure of the image of $\pi_1(k)$ acting on $V_{\ell}$ is connected, then for every integer $d\geq 0$ there exists a constant $D$, depending only on $d$ and $\underline T$,  such that, for every finite field extension $k\subseteq k'$ of degree $\leq d$, one has $(T_{\ell}/\ell)^{\pi_1(k)}=(T_{\ell}/\ell)^{\pi_1(k')}$ for every prime $\ell \geq D$. To prove this, one exploits again independence results, not in the $\ell$-adic setting but in the ultrafilter setting, recently obtained by Cadoret-Hui-Tamagawa in \cite{ultra} and by Cadoret in \cite[Section 15]{Weil2Anna}.

Smooth proper base change theorem, the Weil conjectures (\cite{Weil2}) and the independence of $\ell$ of homological equivalence for divisors show that $\{T_{\ell}(Br(Y_{\overline k})):=\varprojlim_nBr(Y_{\overline k})[\ell^n]\}_{\ell\neq p}$ is a compatible system. In this setting, Propositions \ref{maximalcompatible} and \ref{adeliccompatible} are the positive characteristic analogues of \cite[Propositions 5.4 and 5.5]{skoro}, hence we can conclude the proof of Theorem \ref{mainBrauer} adjusting the arguments in \cite[Section 5.4]{skoro}.
\subsection{Organization of the paper}
In Section \ref{connected} we prove Theorems \ref{maximalcompatible} and \ref{adeliccompatible}. In Section \ref{sectionbrauer} we apply Theorems \ref{maximalcompatible} and \ref{adeliccompatible} to representations coming from geometry and we prove Theorem \ref{mainBrauer}. We end the paper in Section \ref{abelianvarietiesection} discussing applications to abelian varieties. 
\subsection{Acknowledgements}
This paper is part of the author Ph.D. project under the supervision of Anna Cadoret. He thank her for many (many) useful discussions and insights and for her careful re-(re)-reading. The author is also grateful to Akio Tamagawa for suggesting the counterexample in Footnote \ref{footnote}.
\subsection{Conventions and notation}
For the rest of the paper $k$ is a finitely generated field of characteristic $p>0$ with algebraic closure $k\subseteq \overline k$. We write $\F_q$ (resp. $\F$) for the algebraic closure of $\F_p$ in $k$ (resp. $\overline k$) and $k_{\F}:=k\otimes_{\F_q}\F\simeq k\F\subseteq \overline k$. If $R$ is a commutative ring, $A$ an $R$-module and $n,m$ integers $\geq 0$, set
$$T^{n,m}(A):=\underbrace{A\otimes_R \ldots \otimes_R A}_{n\text{ times}}\otimes_R \underbrace{A^{\vee}\otimes_R \ldots \otimes_R A^{\vee}}_{m \text{ times}}.$$
If $G$ is an algebraic group over a field, write $G^0$ for  its neutral component and $\pi_0(G)$ for the group of connected components. Write $\ell_{0}=3$ (resp. $\ell_0=2$) if $p\neq 3$ (resp. $p=3$) and $s_{\ell}=\ell$ (resp. $s_{\ell}=4$) if $\ell\neq 2$ (resp. $\ell=2$).
\section{Forms of representations}\label{connected}
\subsection{Proof of Proposition \texorpdfstring{\ref{maximalcompatible}}-}\label{connected component}
Before proving Proposition \ref{maximalcompatible}, we collect a couple of preliminary propositions.
\subsubsection{A group theoretical proposition}
Let $T$ be a free $\Zl$-module of rank $r$ and let $\Pi\subseteq GL(T)$ be a closed subgroup. Write $V:=T\otimes \Ql$ and let $G\subseteq GL(V)$ be the Zariski closure of $\Pi$. Then:
\begin{proposition}\label{connectedcomponent}
$|\pi_0(G)|\leq |GL_r(\Z/s_{\ell})|$
\end{proposition}
\proof
Write $G^{red}$ for the Zariski closure of the image of $\Pi$ acting on the $\Pi$-semisimplification of $V$. Since the kernel of the natural surjection $G\rightarrow G^{red}$ is unipotent hence connected, it induces an isomorphism $\pi_0(G)\simeq \pi_0(G^{red})$. So, one may assume that $G$ is reductive. Write $H:=Ker(\Pi\rightarrow GL(T/s_{\ell}))$. Since $[\Pi:H]\leq |GL_r(\Z/s_{\ell})|$ and $H$ acts trivially on $GL(T/s_{\ell})$, Lemma \ref{connectedmodl} below concludes the proof. \endproof
\begin{lemma}\label{connectedmodl}
If $G$ is reductive and the action of $\Pi$ on $T/s_{\ell}$ is trivial, then $G$ is connected.
\end{lemma}
\proof
By \cite[Lemma 2.3]{larsenpink}, it is enough to show that, for every irreducible representation $W$ of $GL(V)$ one has $W^{G}=W^{G^0}$. 
Since $GL(V)$ is reductive, by \cite[Proposition 3.1]{delignemilne} every irreducible representation of $GL(V)$ is a sub module of $T^{n,m}(V)$ and hence it is enough to show that for every integers $n,m\geq 0$ $$T^{n,m}(V)^{G}=T^{n,m}(V)^{G^0}.$$ 
The $\Zl$-module $T^{n,m}(T)$ is a $\Pi$-invariant $\Zl$-lattice in $T^{n,m}(V)$ and $\Pi$ acts trivially on $T^{n,m}(T)/s_{\ell}=T^{n,m}(T/s_{\ell})$, so that, by \cite[Lemma 2.1]{ghost}, for every open subgroup $U\subseteq \Pi$ one has $$Hom_{\Pi}(\Ql,T^{n,m}(V))=Hom_{U}(\Ql,T^{n,m}(V)).$$
Applying this to $U:=Ker(\Pi\twoheadrightarrow \pi_0(G))$, one gets
$$T^{n,m}(V)^{G}=Hom_{\Pi}(\Ql,T^{n,m}(V))=Hom_{U}(\Ql,T^{n,m}(V))=T^{n,m}(V)^{G^0}.$$
\endproof
\subsubsection{Independence}
Let $?\in \{\emptyset,\F\}$. Let $\underline T$ be a $\pi_1(k)$-compatible system of finitely generated $\Z_{\ell}$-modules of rank $r$ and write $G_{\ell,?}$ for the Zariski closure of the image of $\pi_1(k_?)$ acting on $V_{\ell}:=T_{\ell}\otimes \Ql$.
\begin{corollary}\label{mainconnected}
For every prime $\ell\neq p$ one has $\pi_0(G_{\ell,?})|\leq |GL_{r}(\Z/\ell_0)|$.
\end{corollary}
\proof
By Lemma \ref{connectedcomponent}, it is enough to show that if $G_{\ell_0,?}$ is connected then $G_{\ell,?}$ is connected for every prime $\ell\neq \ell_0$. By definition of a compatible system, there exists a smooth geometrically connected $\F_q$-variety $\mathcal K$ with generic point $Spec(k)\rightarrow \mathcal K$ such that, for every prime $\ell\neq p$, the action of $\pi_1(k)$ on $T_{\ell}$ factors trough the surjection $\pi_1(k)\twoheadrightarrow \pi_1(\mathcal K)$. So it is enough to show the corresponding statement for the actions of $\pi_1(\mathcal K)$ and $\pi_1(\mathcal K_{\F})$ on $V_{\ell}$. This follows from Fact \ref{connectedisind} below.\endproof
\begin{fact}\label{connectedisind}
$G_{\ell_0,?}$ is connected if and only if $G_{\ell,?}$ is connected.
\end{fact}
\proof
To prove Fact \ref{connectedisind} one can replace $V_{\ell}$ with its $\pi_1(\mathcal K)$-semisimplification. So we may and do assume that $V_{\ell}$ is semisimple as $\pi_1(\mathcal K)$-module, hence as $\pi_1(\mathcal K_{\F})$-module.
Then, arguing as in Lemma \ref{connectedmodl}, it is enough to show that for every integers $n,m\geq 0$ one has 
$$T^{n,m}(V_{\ell})^{G_{\ell, ?}}=T^{n,m}(V_{\ell_0})^{G_{\ell_0,?}}.$$
By \cite{lafforgue} and \cite{drinfeld} every semisimple $\pi_1(\mathcal K)$-modules is direct sum of its pure components (see \cite[Theorem 3.5.5]{marcuzzo} for more details) so that one reduces to the situation in which $V_{\ell_0}$ and $V_{\ell}$ are pure. Then, by the theory of weights (\cite{Weil2}), the dimensions of $T^{n,m}(V_{\ell})^{G_{\ell, ?}}$ and $T^{n,m}(V_{\ell_0})^{G_{\ell_0,?}}$, can be read on the L-functions of $T^{n,m}(V_{\ell})^{\vee}(d)$ and $T^{n,m}(V_{\ell_0})^{\vee}(d)$, where $d$ is the dimension of $\mathcal K$ (see \cite[Proposition 3.4.11]{marcuzzo} for more details). Since $T^{n,m}(V_{\ell})$ and $T^{n,m}(V_{\ell_0})$  are compatible, this concludes the proof.
\endproof
 \begin{remark}
 Fact \ref{connectedisind} is proved in \cite{serreribet} if $?=\emptyset$ and in \cite[Theorem 2.2]{larsenpink} if $?=\F$ and $V_{\ell}$ is pure.
 \end{remark}
\subsubsection{Proof of Theorem \texorpdfstring{\ref{maximalcompatible}}- }\label{proofmaximal}
Keep the notation as in the statement of Proposition \ref{maximalcompatible} and fix $?\in \{\emptyset, \F\}$. 
We can replace $T_{\ell}$ with $T_{\ell}/(T_{\ell})_{tors}$ and $U_{\ell}$ with $U_{\ell}/(U_{\ell})_{tors}$, hence assume that $T_{\ell}$ and $U_{\ell}$ are torsion free. Since $\underline T$ and $\underline U$ are compatible systems, $\{H_{\ell}:=T^{\vee}_{\ell}\otimes U_{\ell}\}_{\ell\neq p}$ is a compatible system as well. By Corollary \ref{mainconnected}, there exists a finite field extension $k_?\subseteq k_{\underline U}$ of degree $\leq |GL_{r^2}(\Z/s_{\ell_0})|$ such that the Zariski closure $G_{\ell}$ of the image of $\pi_1(k_{\underline U})$ acting on $H_{\ell}\otimes \Ql$ is connected for every prime $\ell\neq p$. We claim that $k_{\underline U}$ satisfies the conclusion of Proposition \ref{maximalcompatible}. By assumption, there exists a finite extension $k\subseteq k_{\ell}$ and an isomorphism $$\psi_{\ell}\in H_{\ell}^{\pi_1(k_{\ell})}\subseteq H_{\ell}^{\pi_1(k_{\ell}k_{\underline U})},$$
hence it is enough to show that $H_{\ell}^{\pi_1(k_{\underline U})}=H_{\ell}^{\pi_1(k_{\ell}k_{\underline U})}.$
Since $H_{\ell}^{\pi_1(k_{\ell}k_{\underline U})}/H_{\ell}^{\pi_1(k_{\underline U})}$ is torsion free, it is enough to show that $(H_{\ell}\otimes \Ql)^{\pi_1(k_{\underline U})}=(H_{\ell}\otimes \Ql)^{\pi_1(k_{\ell}k_{\underline U})}$
and this follows from the facts that $k_{\underline U}\subseteq k_{\ell}k_{\underline U}$ is a finite field extension and $G_{\ell}$ is connected. This concludes the proof.
\subsection{Proof of Proposition \texorpdfstring{\ref{adeliccompatible}}-}\label{adelic}
Keep the notation as in the statement of Proposition \ref{maximalcompatible} and fix $?\in \{\emptyset, \F\}$. Write
$$V_{\ell}:=T_{\ell}\otimes \Ql;\quad M_{\ell}:=T_{\ell}\otimes \Ql/\Zl; \quad \overline T_{\ell}:=T_{\ell}/\ell.$$
\subsubsection{Preliminary reduction}\label{mprem}
Write $G_{\ell,?}$ for the Zariski closure of the image $\Pi_{\ell,?}$ acting on $V_{\ell}$.
By Corollary \ref{mainconnected} and replacing $k_?$ with a finite field extension of degree $\leq |GL_{r}(\Z/s_{\ell_0})|$, one may assume that $G_{\ell,?}$ is connected for every prime $\ell\neq p$. Since by assumption there are at most finitely many $T_{\ell}$ with torsion and these are finitely generated $\Z_{\ell}$-modules, we may replace $T_{\ell}$ with $T_{\ell}/(T_{\ell})_{tors}$  hence assume that $T_{\ell}$ is torsion free for every prime $\ell\neq p$. The proof of Proposition \ref{adeliccompatible} is the combination of the following two claims and the arguments in Section \ref{end}.
\newline
\textbf{Claim 1:} For every integer $d\geq 1$ and for every prime $\ell\neq p$, there exists a constant $A_{\ell}:=A(d,\ell,\underline T)$ such that, for every finite field extension $k\subseteq k'$ of degree $\leq d$, one has $[M_{\ell}^{\pi_1(k')}:M_{\ell}^{\pi_1(k_?)}]\leq A_{\ell}$.
\newline
\textbf{Claim 2:} For every integer $d\geq 1$, there exists a constant $D:=D(\underline T,d)$ such that, for every prime $\ell\geq D$ and for every finite field extension $k\subseteq k'$ of degree $\leq d$, one has $\overline T_{\ell}^{\pi_1(k')}=\overline T_{\ell}^{\pi_1(k_?)}$.
\subsubsection{Proof of Claim 1}\label{ml}
Since $\Pi_{\ell,?}$ is a compact $\ell$-adic Lie group, it is topologically finitely generated and hence it has finitely many open subgroups of bounded index. So it is enough to show that if $U\subseteq \Pi_{\ell,?}$ is an open subgroup then 
$[M_{\ell}^{U}:M_{\ell}^{\Pi_{\ell,?}}]<+\infty$.
This follows from \cite[Lemma 3.3.2]{brauer} and the connectedness of $G_{\ell,?}$. To the reader convenience, we briefly recall the argument.

Since $G_{\ell,?}$ is connected, one has $V_{\ell}^{\Pi_{\ell,?}}=V_{\ell}^{U}$ and $T_{\ell}^{\Pi_{\ell,?}}=T_{\ell}^{U}.$ The exact sequence 
$$0\rightarrow T_{\ell}\rightarrow V_{\ell}\rightarrow M_{\ell}\rightarrow 0$$
induces a commutative diagram with exact rows:
\begin{center}
\begin{tikzcd}
0\arrow{r}& V_{\ell}^{\Pi_{\ell,?}}/T_{\ell}^{\Pi_{\ell,?}}\arrow[equal]{d}\arrow{r} & M_{\ell}^{\Pi_{\ell,?}}\arrow{r}\arrow[hook]{d} & H^1(\Pi_{\ell,?},T_{\ell})\arrow{d} \\
0\arrow{r}& V_{\ell}^{U}/T_{\ell}^{U}\arrow{r} & M_{\ell}^{U}\arrow{r}{\Delta} & H^1(U,T_{\ell}) 
\end{tikzcd}
\end{center}
So $M_{\ell}^{U}/M_{\ell}^{\Pi_{\ell,?}}$ is a quotient of the image of $\Delta$. But $\Delta$ has finite image since $M_{\ell}^{U}$ is torsion and $H^1(U,T_{\ell})$ is a finitely generated $\Z_{\ell}$-module by \cite[Proposition 9]{serre}.
\subsubsection{Proof of Claim 2}
For any finite field extension $k_?\subseteq k'$, consider the images $\Pi_{k'}\subseteq \Pi_{?}$ of $\pi_1(k')\subseteq \pi_1(k_?)$ acting on $\overline T:=\prod_{\ell\neq p}\overline T_{\ell}$.
By definition of a compatible system, there exists a smooth geometrically connected $\F_q$-variety $\mathcal K$ with generic point $Spec(k)\twoheadrightarrow \mathcal K$ such that, for every prime $\ell\neq p$, the action of $\pi_1(k)$ on $T_{\ell}$ factors trough the canonical surjection $\pi_1(k)\rightarrow \pi_1(\mathcal K)$. By the Grothendieck-Ogg-Shafarevich formula, there exists a connected \etale cover $\mathcal K'\rightarrow \mathcal K$ such that the action of $\pi_1(\mathcal K')\subseteq \pi_1(\mathcal K)$ on $\overline T$ factors trough the tame fundamental group of $\mathcal K'$; see the proof of \cite[Lemma 12.3.1]{Weil2Anna}. Since the tame fundamental groups of $\mathcal K'$ and of every connected component of $\mathcal K'_{\F}$ are topologically finitely generated, this implies that $\Pi_{?}$ is topologically finitely generated. Hence the group $\Pi_{?}$ has finitely many open subgroups of index $\leq d$. So there are only finitely many possibilities for the inclusions $\Pi_{k'}\subseteq \Pi_{?}$, while $k_?\subseteq k'$ is varying among the finite field extensions of degree $\leq d$. So, to prove Claim 2, it is enough to show\footnote{This is not a formal consequence of $[\pi_1(k_{?}):\pi_1(k')]$ being finite, as the example $\{1\}\subseteq \{1,-1\}\subseteq \prod_{\ell\neq p}GL(\overline T_{\ell})$ shows.} that, for every finite field extension $k_{?}\subseteq k'$ of degree $\leq d$, there exists a constant $D':=D'(\underline T,k')$ such that for $\ell\geq D'$ one has $\overline T_{\ell}^{\pi_1(k')}=\overline T_{\ell}^{\pi_1(k_?)}$.

Let $\mathcal L$ be the set of prime $\neq p$ and write $F:=\prod_{\ell\in \mathcal L}\F_{\ell}$. We use the formalism of ultrafilters\footnote{In this note an ultrafilter will always mean a non-principal ultrafilter.} on $\mathcal L$; see \cite[Appendix]{ultra}. To every ultrafilter $\mathfrak u$ on $\mathcal L$ one associates a maximal ideal $\mathfrak m_{\mathfrak u}$ of $F$ and writes $F_{\mathfrak u}:=F/\mathfrak m_{\mathfrak u}$ for the characteristic zero residue field. The actions of $\pi_1(k_{?})$ and $\pi_1(k')$ on $\overline T$ induces actions on $T_{\mathfrak u}:=T\otimes_F F_{\mathfrak u}$.
Since $\Pi_{?}$ and $\Pi_{k'}$ are topologically finitely generated groups, by \cite[Lemma 4.3.3]{ultra} and \cite[Lemma 4.4.2]{ultra} it is enough to show that $T_{\mathfrak u}^{\pi_1(k_{?})}=T_{\mathfrak u}^{\pi_1(k')}$
for every ultrafilter $\mathfrak u$. Write $G_{\mathfrak u,?}$ and $G_{\mathfrak u,k'}$ for the Zariski closures of the images of $\pi_1(k_?)$ and $\pi_1(k')$ acting on $T_{\mathfrak u}$. Since $T_{\mathfrak u}^{\pi_1(k_{?})}=T_{\mathfrak u}^{G_{\mathfrak u,?}}$ and $T_{\mathfrak u}^{G_{\mathfrak u,k'}}=T_{\mathfrak u}^{\pi_1(k')}$, it is enough to show that the natural inclusion $G_{\mathfrak u,k'}\subseteq G_{\mathfrak u,?}$ is an equality. Since $\pi_1(k')\subseteq \pi_1(k_?)$ has finite index, one has $G_{\mathfrak u,k'}^0=G_{\mathfrak u,?}^0$ hence it is enough to show that $G_{\mathfrak u,?}$ is connected. This follows from the fact that $G_{\ell}$ is connected by preliminary reduction and Fact \ref{connectedultrafilter} below.
\begin{fact}\label{connectedultrafilter}
The group $G_{\ell,?}$ is connected if and only if $G_{\mathfrak u,?}$ is connected.
\end{fact} 
\proof
If $?=\emptyset$ this is proved in \cite[Theorem 1.3.1]{ultra} and if $?=\F$ this is proved in \cite[Corollary 15.1.2]{Weil2Anna}. \endproof
 \numberwithin{equation}{subsubsection}
\subsubsection{End of the proof}\label{end}
To conclude the proof of Proposition \ref{adeliccompatible}, fix a finite field extension $k_?\subseteq k'$ of degree $\leq d$. Up to replacing $d$ with $d!$ we may restrict to finite Galois extensions $k_?\subseteq k'$, so that $\pi_1(k')\subseteq \pi_1(k_?)$ is a normal subgroup. By Claim 1, it is enough to show that there exists a constant $A:=A(\underline T,d)$ such that for $\ell\geq A$ one has $M_{\ell}^{\pi_1(k_{?})}=M_{\ell}^{\pi_1(k')}$
and, by Claim 2, there exists a constant $D:=D(\underline T,d)$ such that for $\ell\geq D$ one has
$\overline T_{\ell}^{\pi_1(k_{?})}=\overline T_{\ell}^{\pi_1(k')}$. We claim that $A:=\max(D,d+1)$ has the desired property.

Since $M_{\ell}=\varinjlim_nM_{\ell}[\ell^n]$, it is enough to show that for $\ell\geq A$ and every $n\geq 1$ one has $M_{\ell}[\ell^{n}]^{\pi_1(k_{?})}=M_{\ell}[\ell^{n}]^{\pi_1(k')}$. 
For this, one argues by induction on $n$, the case $n=1$ being the definition of $D$.
For $n>1$, since $T_{\ell}$ is torsion free, there is a $\pi_1(k_{?})$-invariant identification $M_{\ell}[\ell^{n}]\simeq T_{\ell}/\ell^n$ and a $\pi_1(k_{?})$-equivariant exact sequence
$$0\rightarrow \overline T_{\ell}\rightarrow T_{\ell}/\ell^n\rightarrow T_{\ell}/\ell^{n-1}\rightarrow 0.$$
Combined with the inflation-restriction exact sequence for the normal inclusion $\pi_1(k')\subseteq \pi_1(k_{?})$, this induces a commutative exact diagram
\begin{center}
\begin{tikzcd}
& &&& H^1(\pi_1(k_{?})/\pi_1(k'),\overline T_{\ell}^{\pi_1(k')})\arrow{d}\\\
0\arrow{r}& \overline T_{\ell}^{\pi_1(k_{?})}\arrow{r}\arrow{d}{\simeq}&(T_{\ell}/\ell^{n})^{\pi_1(k_{?})}\arrow{r}\arrow[hook]{d} & \arrow{r}\arrow{d}{\simeq} (T_{\ell}/\ell^{n-1})^{\pi_1(k_{?})}& H^1(\pi_1(k_{?}),\overline T_{\ell})\arrow{d}\\
0\arrow{r}& \overline T_{\ell}^{\pi_1(k')}\arrow{r}&(T_{\ell}/\ell^{n})^{\pi_1(k')}\arrow{r} & \arrow{r} (T_{\ell}/\ell^{n-1})^{\pi_1(k')}& H^1(\pi_1(k'),\overline T_{\ell}).
\end{tikzcd}
\end{center}
By the induction hypothesis the first and the third vertical arrows are isomorphisms for $\ell\geq A$. By elementary diagram chasing it is enough to show that $H^1(\pi_1(k_{?})/\pi_1(k'),\overline T_{\ell}^{\pi_1(k')})=0$. But since $\overline T_{\ell}^{\pi_1(k_{?})}=\overline T_{\ell}^{\pi_1(k')}$ one has 
$$H^1(\pi_1(k_{?})/\pi_1(k'),\overline T_{\ell}^{\pi_1(k')})=H^1(\pi_1(k_{?})/\pi_1(k'),\overline T_{\ell}^{\pi_1(k_{?})})=Hom(\pi_1(k_{?})/\pi_1(k'),(\Z/\ell)^r)=0$$
where the last equality follows from the fact that $\ell>d=|\pi_1(k_{?})/\pi_1(k')|$.\endproof
\section{Proof of Theorem \texorpdfstring{\ref{mainBrauer}}- }\label{sectionbrauer}
\subsection{Proof of Theorem \texorpdfstring{\ref{mainBrauer}}- }\label{proof brauer}
Retain the notation and the assumption of Proposition \ref{mainBrauer}. For every finite field extension $k\subseteq k'$ and every $(\overline k/k')$-form $Y$ of $X$, write
$Y_{\overline k}:=Y\times_{k'}\overline k$ and 
\begin{gather*}
T_{\ell}(Y):=\varprojlim_nBr(Y_{\overline k})[\ell^n]; \quad M_{\ell}(Y):=T_{\ell}(Y)\otimes \Ql/\Zl;  \quad \underline{M}(Y):=\prod_{\ell\neq p} M_{\ell}(Br(Y_{\overline{k}}));\\
H^2_{\ell}(Y):=H^2(Y_{\overline k},\Zl(1));\quad \underline{H}^i(Y):=\{H^2_{\ell}(Y)\}.
\end{gather*}
\subsubsection{Reducing to the Tate module of the Brauer group}
Recall (see e.g. the proof of \cite[Proposition 2.1.1]{brauer}) that there is a $\pi_1(k')$-equivariant exact sequence
$$ 0\rightarrow M_{\ell}(Br(Y_{\overline{k}}))\rightarrow Br(Y_{\overline{k}})[\ell^{\infty}]\rightarrow H^3(Y_{\overline k},\Zl(1))[\ell^{\infty}]\rightarrow 0.$$
Since
\begin{itemize}
\item for every prime $\ell\neq p$, the group $H^3(Y_{\overline k},\Zl(1))[\ell^{\infty}]=H^3(X_{\overline k},\Zl(1))[\ell^{\infty}]$ is finite (of cardinality depending only on X) and
\item for $\ell\gg 0$ (depending only on $X$) one has $H^3(Y_{\overline k},\Zl(1))[\ell^{\infty}]=H^3(X_{\overline k},\Zl(1))[\ell^{\infty}]=0 $ (\cite{gabber});
\end{itemize}
it is enough to prove Theorem \ref{mainBrauer} replacing $Br(Y_{\overline k})[p']$ with $\underline{M}(Y).$
\subsubsection{Compatibility}\label{sectionmotivicbrauer}
We now prove that $\underline T(Y)$ is a compatible system of $\pi_1(k')$-modules.  Write $NS(Y_{\overline k})$ for the N\'eron-Severi group of $Y_{\overline k}$. By the Kummer exact sequence
$$ 0\rightarrow NS(Y_{\overline k})\otimes \Zl \rightarrow H^2_{\ell}(Y)\rightarrow T_{\ell}(Y)\rightarrow 0,$$
it is enough to show that $\underline H^2$ and $\underline{NS}(Y):=\{NS(Y_{\overline k})\otimes \Zl\}_{\ell\neq p}$  are compatible systems of $\pi_1(k')$-modules.
Write $\F_{q'}$ for the algebraic closure of $\F_q$ in $k'$.
By spreading out, there exists a geometrically connected smooth $\F_{q'}$-variety $\mathcal K'$, with generic point $\eta':Spec(k')\rightarrow \mathcal K'$, and a smooth proper morphism $\mathfrak f:\mathcal Y\rightarrow \mathcal K'$ fitting into a commutative cartesian diagram:
 \begin{center}
\begin{tikzcd}
Y\arrow{r}\arrow{d}\arrow[phantom]{rd}{\Box} & \mathcal Y\arrow{d}{\mathfrak f}\\
Spec(k')\arrow{r}{\eta'} & \mathcal K'.
\end{tikzcd}
\end{center}
By smooth proper base change, the action of $\pi_1(k')$ on $H^2_{\ell}(Y)$ factors trough the surjection $\pi_1(k')\twoheadrightarrow \pi_1(\mathcal K')$ and by \cite{Weil2} the collection $\underline H^{2}(Y)$ is a $\Q$-rational compatible system. Since homological and algebraic equivalences coincide rationally for divisors, $NS(Y_{\overline k})\otimes \Q$ identifies with the image of the cycle class map $c_{Y_{\overline k}}:Pic(Y_{\overline k})\otimes \Q\rightarrow H^2_{\ell}(Y)\otimes \Ql$. So $\underline{NS}(Y)$ is a compatible system of $\pi_1(k')$-modules, hence $\underline T(Y)$ is a compatible system of $\pi_1(k')$-modules as well.
\subsubsection{End of the proof}
So we can apply Propositions \ref{maximalcompatible} and \ref{adeliccompatible} to $\underline T(Y)$. Hence, to conclude the proof, we have just to adjust the arguments in \cite[Section 5.4]{skoro}, replacing \cite[Propositions 5.4 and 5.5]{skoro} with Propositions \ref{maximalcompatible} and \ref{adeliccompatible}. Write $r:=Rank_{\Zl}(T_{\ell}(X))^2$ and set $B_X:=|GL_r(\Z/\ell_0)|$.
By Proposition \ref{maximalcompatible} for $X_{k'}$ there exists a finite field extension $k'\subseteq k_Y$ of degree $\leq B_X$ such that there is an  $\pi_1(k_Y)$-equivariant isomorphism 
$\underline{M}(Y)\simeq \underline{M}(X).$
Then one has:
$$\underline{M}(X)^{\pi_1(k)}\subseteq \underline{M}(X)^{\pi_1(k_Y)}\simeq \underline{M}(Y)^{\pi_1(k_Y)}\supseteq \underline{M}(Y)^{\pi_1(k')}.$$
Since $T(X,\ell)$ holds for every prime $\ell\neq p$, by Fact \ref{finiteness} the group $\underline{M}(X)^{\pi_1(k_Y)}$ is finite. Hence it is enough to show that, for every integer $d\geq 1$, there exists a constant $C:=C(X,d)$ such that for every finite field extension $k\subseteq k''$ of degree $\leq d$ one has $\underline{M}(X)^{\pi_1(k'')}\leq C.$
To prove this, one may replace $k$ with a finite extension and then apply Proposition \ref{adeliccompatible} to conclude. \endproof
\subsection{Further remarks}\label{abelianvarietiesection}
Let $k$ be an infinite finitely generated field of characteristic $p\geq 0$.
\subsubsection{Torsion of abelian varieties}\label{Sectionabelian}
Let $X$ be a $k$-abelian variety of dimension $g$.  By the Lang-N\'eron theorem \cite{LN}, the group $X(k')_{tors}$ is finite for every finite field extension $k\subseteq k'$ and, if $X$ has no isotrivial geometric isogeny factors, then the same is true for every field extension of $k_{\F}$. One can use Theorems \ref{maximalcompatible} and \ref{adeliccompatible} with the techniques in Section \ref{proof brauer} to prove uniform boundedness results for the torsion of the forms of abelian varieties. More precisely, one can prove that for every integer $d\geq 1$, (resp. if $X$ has no isotrivial geometric isogeny factors) there exists an integer $C:=C(X,d)$ such that $|Y(k')|\leq C$ for every finite extension of fields $k\subseteq k'$ (resp. $k_{\F}\subseteq k'$) of degree $\leq d$ and every $k'$-abelian variety $Y$ that is a $(\overline k/k')$ form of $X$. We conclude pointing out that the statement for abelian varieties over $k$ follows also from the Lang-Weil bound and the specialization theory for torsion of abelian varieties. 
\subsubsection{Abelian varieties with CM}\label{sectionCM}
Recall that a $k$-abelian variety $X$ has complex multiplication (or $CM$ for short) if the image of the representation $\pi_1(k)\rightarrow GL(T_{\ell}(X))$ contains an abelian open subgroup. In characteristic zero, Orr-Skorobogatov (\cite[Corollary C.2]{skoro}) prove that there is a constant $C = C(d,g)$ such that $|Br(\Xbar)^{\pi_1(k)}|\leq C$ for every $g$-dimensional abelian variety with CM defined over a number field $k$ of degree $\leq d$. This result is a consequence of the characteristic zero analogue \cite[Theorem 5.1]{skoro} of Theorem \ref{mainBrauer} and of the fact (\cite[Theorem A]{skoro}) that there are only finitely many $\overline {\Q}$-isomorphism classes of $g$-dimensional abelian varieties with $CM$ defined over a number field of degree $\leq d$. Unfortunately, as Akio Tamagawa pointed out to us, the positive characteristic analogue of \cite[Theorem A]{skoro} is false: if $X$ is the product of $g>1$ supersingular elliptic curves, the $k$-isogeny class of $X$ contains infinitely many\footnote{\label{footnote}Indeed, there is an inclusion $\alpha_p^2\subseteq X$. Since $k$ is infinite, the set $I:=Hom_k(\alpha_p,\alpha_p\times \alpha_p)/Aut_k(\alpha_p)\simeq \mathbb P^1(k)$ is infinite.
For each $i\in I$ define $f_i:X\rightarrow X_i:=X/i(\alpha_p)$. Assume by contradiction that the $X_{i,\overline k}$ fall into finitely isomorphism many classes. Then there exist $i_0$ and an infinite subset $J\subseteq I$ such that, for every $j\in J$, there is an isomorphism $g_j:X_{j,\overline k}\rightarrow X_{i_0,\overline k}$. Then, $g_{j}\circ f_{j}:X_{\overline k}\rightarrow X_{i_0,\overline k}$ is a map of degree $p$. Since there are only finitely many maps $X_{\overline k}\rightarrow X_{i_0,\overline k}$ of degree $p$ , there exists an infinite subset $J'\subseteq J$ such that  $g_{j}\circ f_j=g_{j'}\circ f_{j'}$ for every $j,j'\in J'$. But this implies $j(\alpha_p)=j'(\alpha_p)$ and this is a contradiction.} $k$-abelian varieties that are not isomorphic over $\overline k$. So there is no hope to deduce directly from Theorem \ref{mainBrauer} the analogue of \cite[Corollary C.2]{skoro} in positive characteristic. However, a positive characteristic version of \cite[Corollary C.2]{skoro}, via different techniques, has been announced by Marco D'Addezio.

\end{document}